\documentclass[12pt]{article}

\usepackage[a4paper,hmargin=2cm,vmargin=3cm]{geometry}
\usepackage{fancyhdr,graphicx,lastpage}
\usepackage{listings}
\usepackage{amsmath}
\usepackage{amssymb}
\usepackage{amsthm}
\usepackage{enumerate}
\usepackage{parskip}
\usepackage{multirow}
\usepackage{url}
\usepackage{makeidx}
\usepackage{mathrsfs}
\makeindex
\newtheoremstyle{krisjan}
{}
{}
{}
{}
{\bf}
{.}
{3pt}
{}

\theoremstyle{krisjan}
\newtheorem{thm}{Theorem}
\newtheorem{prop}[thm]{Proposition}
\newtheorem{lem}[thm]{Lemma}

\title{A note on Products of Nilpotent Matrices}
\author{Christiaan Hattingh\footnote{University of South Africa, Department of Mathematics}
\footnote{e-mail address: krisjan@simulasie.co.za}}
\begin{document}
\maketitle

\begin{abstract}
Let $\mathscr{F}$ be a field. A matrix $A \in M_n(\mathscr{F})$ is a product of two nilpotent matrices if and only if it is singular, except if $A$ is a nonzero nilpotent matrix of order $2 \times 2$. This result was proved independently by Sourour \cite{sour}  and Laffey \cite{laffey}. While these results remain true and the general strategies and principles of the proofs correct, there are certain problematic details in the original proofs which are resolved in this article. A detailed and rigorous proof of the result based on Laffey's original proof \cite{laffey} is provided, and a problematic detail in the Proposition which is the main device underpinning the proof by Sourour is resolved.
\end{abstract}
\section{Introduction}
First I will fix some notation. A standard basis vector, which is a column vector with a one in position $i$ and zeros elsewhere, is indicated as $e_i$. A matrix with a one in entry $(i,j)$ and zeros elsewhere is indicated as $E_{(i,j)}$. Scripted letters such as $\mathscr{F}$ indicate a field, $M_n(\mathscr{F})$ indicate the set of all matrices of order $n \times n$ over the field $\mathscr{F}$. A block diagonal matrix with diagonal blocks $A$ and $B$ (in that order from left top to bottom right) is indicated as $\text{Dg}[A,B]$. The notation $[a,b,\ldots,k]$ indicates a matrix with the column vectors $a, b, \ldots, k$ as columns. A simple Jordan block matrix of order $k$, and with eigenvalue $\lambda$ is indicated as $J_k(\lambda)$. In this text the convention of ones on the subdiagonal is assumed for the Jordan canonical form, unless otherwise indicated.

\section{Problematic details in the original proofs}
The proof by Laffey relies on preceding results by Wu \cite{wu}. It was mentioned before \cite{bukovsek} that the right factor on p.229, the last factorization \cite[Lemma 3]{wu}, is in fact not nilpotent for certain values of $k$. Explicitly, the given factorization of $\text{Dg}[J_k(0),J_2(0)]$ is invalid for odd $k$, since the matrix 
\[ \begin{bmatrix} \mathbf{0} & J_2(0) \\ J_k(0) & \mathbf{0} \end{bmatrix} \] is not nilpotent when $k=7$.  

Now in the original proof by Laffey \cite{laffey} there is also an error in the factorization given at the top of page 99:
\[ \left [ {\renewcommand{\arraystretch}{1.2}\begin{array}{c|c} 
N_1 & \mathbf{0} \\
\hline
\begin{matrix} 
0 & \cdots & 0 & l_{11} \\
0 &  & \vdots &  l_{21} \\ 
\vdots &  & & \vdots \\
0 & \cdots & 0 & l_{k1} 
\end{matrix} & \begin{matrix} 
0 & 0 & \cdots & 0 \\
l_{22} & 0 &  &  0 \\ 
\vdots & \ddots & \ddots & \vdots \\
l_{k2} & \cdots & l_{kk} & 0 
\end{matrix}
\end{array}} \right ]
\left [ {\renewcommand{\arraystretch}{1.2}\begin{array}{c|c} 
N_2 &\begin{matrix} 
0 & \cdots &   & 0  \\
\vdots &  &    & \vdots \\ 
0 & \cdots  &  & 0 \\
u_{11} & u_{12} &  \cdots & u_{1k} 
\end{matrix}  \\
\hline
\mathbf{0} &
 \begin{matrix}
0  & u_{22} &  \cdots & u_{2k} \\
\vdots & \ddots & \ddots    &  \vdots \\ 
 \quad \, &  & 0 & u_{kk} \\
0 & \cdots & 0 & 0 
\end{matrix}
\end{array}} \right ].\]
If the last column of $N_1$ is not zero (which is possible according to the original proof), then multiplication of the factors result in a block matrix \[ \begin{bmatrix} A_0 & C \\ B & A_1 \end{bmatrix} \] where the matrix $C$ is not necessarily the zero matrix (as indicated in the original proof). Therefore we cannot apply Roth's theorem to prove similarity of the factorization with the original matrix $A$, as claimed.

Lastly I point out an implicit assumption in the proof by Sourour \cite[p.304]{sour}. The main device underpinning the proof is a Proposition, the statement of which is given in Proposition \ref{prop:sourInductSingular} of the main results below.  On p.305 \cite{sour} the statement is made that $\text{N}(PAP) = \mathcal{M}_1 \oplus \text{N}(A)$. However, in the given construction it is assumed that the vector $e_0$ is not in $\mathcal{M}_2$. Suppose $e_0$ is in $\mathcal{M}_2 = \text{R}(P)$: then we will have $N(PAP) = \mathcal{M}_1 \oplus \text{N}(A) \oplus \text{span}\{e_0\}$, as shown in the following example. Choose $A=J_3(0)$ (which is not square-zero) together with $e_0=(1,0,0,0)^T$ and then select $e_0$ as one of the basis vectors of $\mathcal{M}_2$. The form as specified in the statement of the proposition is then not achieved, since $\text{r}(PAP) < \text{r}(AP) = \text{r}(PA)$. In Proposition \ref{prop:sourInductSingular} of the main results below I provide a proof that explicitly ensures $e_0$ is not in $\text{R}(P)$ ($x_0$ in Proposition 1).   

\section{Main results}
I present results which repair the shortcomings listed in the previous section. First I will address the sufficient part of Sourour's Proposition \cite{sour}. Note that Sourour only makes use of the sufficient part when proving the main theorem on nilpotent factorization. Then I present a complete proof of the main result, following Laffey's \cite{laffey} original proof with the necessary corrections.
\subsection{Sourour's Proposition}
\begin{quote}
\begin{prop}[Proposition in \cite{sour}]
Let $\mathscr{F}$ be an arbitrary field. If $A \in M_n(\mathscr{F})$ is not square-zero, then it is similar to a matrix of the form 
\[\begin{bmatrix} \lambda & c^T \\ b & D \end{bmatrix}\] with $\lambda \in \mathscr{F}$, $D \in M_{(n-1)}(\mathscr{F})$, $\text{r}(D) = \text{r}(A) - 1$, $b \in \text{R}(D)$ and $c \in \text{R}(D^T)$.
\label{prop:sourInductSingular}
\end{prop}\end{quote}
{\bf Proof.} First suppose that $A$ is a scalar matrix. In this case the result is immediate, as $A$ is already in the desired form.

Suppose that $A$ is not square-zero, and $A$ is not scalar. Then we can find a vector $x_0$ such that $Ax_0$ and $x_0$ are linearly independent, and $A^2x_0 \neq 0$. Explicitly, if $A$ is full rank, this result is easy to see, as the fact that $A$ is not scalar ensures that there exists a vector $x_0$ such that $Ax_0$ and $x_0$ are linearly independent, and then $A^2x_0 \neq 0$ since $A$ is full rank. Suppose $A$ is not full rank. Since $A$ is not square-zero there exists $x_1 \in \text{R}(A)$ such that $x_1 \notin \text{N}(A)$. Let $Ax_0=x_1$. Now if $x_0$ and $Ax_0 = x_1$ are linearly independent, the claim holds. If this is not the case then $\lambda x_0 = Ax_0$. Replace $x_0$ with $x_0 + n$ where $n$ is a nonzero vector in the null space of $A$: then $A(x_0 + n) = Ax_0 = \lambda x_0$ and since $x_0$ and $x_0 + n$ are linearly independent the claim holds.

Having established the existence of a suitable vector $x_0$, construct a basis for $\mathscr{F}^n$ in the following way 
\[ \alpha_1 = \{x_1, x_1-x_0, n_1, n_2, \ldots, n_k, u_1, u_2, \ldots, u_{(n-k-2)}\},\] where $n_1, n_2, \ldots n_k$ are basis vectors for the null space of $A$, and $u_1, u_2, \ldots, u_{(n-k-2)}$ are arbitrary subject to $\alpha_1$ being a linearly independent set. Let \[\mathcal{V}_1 = \text{span}\{x_1\} \text{ and}\] \[ \mathcal{V}_2 = \text{span}\{x_1-x_0,n_1, n_2, \ldots, n_k,u_1,u_2,\ldots,u_{(n-k-2)}\},\] so that $\mathscr{F}^n = \mathcal{V}_1 \oplus \mathcal{V}_2$. Let $P$ be the projection along $\mathcal{V}_1$ onto $\mathcal{V}_2$. Now the following two results hold:
\[ \text{N}(PA) = \text{span}\{x_0\} \oplus \text{N}(A), \]
\[ \text{N}(AP) = \mathcal{V}_1 \oplus \text{N}(A).\]
The results hold, since by virtue of the construction of $\alpha_1$ and $P$ we have that the null space of $P$ is contained in the range of $A$, and the null space of $A$ is contained in the range of $P$.

Now the key to obtaining the required result is that we must have
\[\text{N}(PAP) = \text{N}(AP) =  \mathcal{V}_1 \oplus \text{N}(A).\]
For this result to be true, we must have $x_1 \notin \text{R}(AP)$. I will now prove this fact. By virtue of the construction of $\alpha_1$ we have $x_0 \notin \text{R}(P)$, for if it was then $x_0 + (x_1-x_0) = x_1 \in \text{R}(P) = \mathcal{V}_2$ which is a contradiction. Since $x_0 \notin \text{R}(P)$ it follows that $Ax_0 = x_1 \notin \text{R}(AP)$, as required.

Applying the rank-nullity theorem to each of the products above we have 
\begin{equation} \text{r}(PA) = \text{r}(AP) = \text{r}(PAP) = \text{r}(A) - 1. \label{eq:sourRank} \end{equation}

Now relative to the basis $\alpha_1$ the matrix representations of the products above are
\[PA \approx \begin{bmatrix} 0 & 0 \\ b & D \end{bmatrix},\]
\[AP \approx \begin{bmatrix} 0 & c^T \\ 0 & D \end{bmatrix},\]
\[PAP \approx \begin{bmatrix} 0 & 0 \\ 0 & D \end{bmatrix}.\]

Combining this result with \eqref{eq:sourRank}, it is immediately apparent that $b$ is in the column space of $D$ and $c$ is in the range of $D^T$, as desired. \hfill $\square$

\subsection{Complete proof following Laffey}
I will now present a rigorous proof following the result by Laffey \cite{laffey}. First let us address the problematic factorization by Wu \cite{wu}. 
\begin{quote}
\begin{lem}[Lemma 3 \cite{wu}]
Let $\mathscr{F}$ be an arbitrary field. Let $k$ be a positive, odd integer. Then the matrix $\text{Dg}[J_k(0),J_2(0)]$ is the product of two nilpotent matrices, each with rank equal to the rank of $\text{Dg}[J_k(0),J_2(0)]$.
\label{lem:j0nilpotent}
\end{lem}
\end{quote}
{\bf Proof.} First consider the case $k=1$, then 
\begin{equation} \text{Dg}[J_1(0),J_2(0)] = E_{(3,1)}E_{(1,2)}, \label{eq:nilpotent_jkj2_1} \end{equation} and it is immediately apparent that both factors are nilpotent and rank 1 as required.

Suppose $k \geq 3$. Then \begin{equation} \text{Dg}[J_k(0),J_2(0)] = \left [ {\renewcommand{\arraystretch}{1.2}\begin{array}{c|c} 
\mathbf{0} &  A_1 \\
\hline
J_2(0) & \mathbf{0} \end{array}} \right ] 
\left [ {\renewcommand{\arraystretch}{1.2}\begin{array}{c|c} 
\mathbf{0} &  \begin{matrix} 1 & 0 \\ 0 & 0 \end{matrix} \\ 
\hline
 A_2 & \mathbf{0} \end{array}} \right ]=N_1N_2, \label{eq:nilpotent_jkj2_odd} \end{equation}
 where 
 \[A_1 = [e_2,\mathbf{0},e_4,e_3,e_6,e_5,\ldots,e_{k-1},e_{k-2},e_k]\] and
 \[A_2 = [e_1,e_4,e_3,e_6,e_5,\ldots,e_{k-1},e_{k-2},e_k,\mathbf{0}].\]
 Notice in particular that for the case $k=3$ we have $A_1=[e_2,\mathbf{0},e_3]$ and $A_2 = [e_1,e_3,\mathbf{0}]$.
 Now the factor $N_1$ has rank $k$ which is the same as the rank of $\text{Dg}[J_k(0),J_2(0)]$. Let 
 \begin{equation}Q_1 = [e_1,e_{k+2},e_{k},e_{k-1},e_{k-4},e_{k-5},e_{k-8}\ldots,e_{3-(-1)^{(k-3)/2}},e_{k+1},e_{k-2},e_{k-3},\ldots,e_{3+(-1)^{(k-3)/2}}]. \label{eq:nilpotent_jkj2_odd_q} \end{equation} Then $Q_1$ is a change-of-basis matrix and 
 \[Q_1^{-1}N_1Q_1 = \text{Dg}[J_{k-2 \cdot (\lfloor (k-3)/4 \rfloor)+(-1)^{(k-3)/2}}(0),J_{k-2 \cdot (1+\lfloor (k-3)/4 \rfloor)}(0)],\] so that $N_1$ is nilpotent.
 
 Finally, it is easy to verify that $N_2$ also has rank $k$. Let 
 \[Q_2 = [e_4,e_8,\ldots,e_{k+1},e_1,e_3, \ldots,e_{k},e_2,e_6, \ldots,e_{k-1},e_{k+2} ] \] when $(k-3)/2$ is even, and 
 \[Q_2 = [e_2,e_6, \ldots, e_{k+1},e_1,e_3, \ldots,e_{k},e_4,e_8,\ldots,e_{k-1},e_{k+2}] \] when $(k-3)/2$ is odd. Then 
 \[Q_2^{-1} N_2 Q_2 = \text{Dg} [J_{k-\lfloor (k-3)/4 \rfloor}(0),J_{2+\lfloor (k-3)/4 \rfloor}(0)],\] confirming that $N_2$ is nilpotent, and completing the proof. \hfill $\square$

To repair the error in Laffey's original proof \cite[p.99]{laffey} we require the following result.
\begin{quote}
\begin{lem}
Let $\mathscr{F}$ be an arbitrary field. Any nilpotent matrix $N \in M_n(\mathscr{F})$, with $n \neq 2$, is similar to the product of two nilpotent matrices, where the first row and last column of the first factor (left factor) is the zero vector, and the last row of the second factor (right factor) is either the zero vector or $e_1^T$ (the vector with a one in the first entry and zeros elsewhere).
\label{lem:col0row0}
\end{lem}
\end{quote}
{\bf Proof.}
Wu \cite[Lemma 3]{wu} present exhaustive configurations of $N$, and factorizations into nilpotent factors for each of these. Now we need to prove that each of these factorizations are of the desired form, possibly using a similarity transformation where needed.

First, I prove that for each configuration of the last block(s) of $N$ we can find a factorization, possibly with a suitable similarity transformation, where the last column of the first factor is the zero vector and the last row of the second factor is either the zero vector or $e_1^T$. 

Suppose $J_k(0)$ with $k \neq 2$ is the last simple Jordan block of $N$ on the diagonal (i.e. this block is not paired with a block $J_2(0)$). If $k=1$ the result is immediate. Suppose $k > 2$: now we need to consider the factorizations \begin{equation} J_k(0) = \left [ {\renewcommand{\arraystretch}{1.2}\begin{array}{c|c|c} 
\mathbf{0} &  \begin{matrix} 0 \\ 0 \end{matrix} & \begin{matrix} 0 \\ 1 \end{matrix}\\
\hline 
I_{k-2} &  \mathbf{0} &  \mathbf{0} \end{array}} \right ]
\left [ {\renewcommand{\arraystretch}{1.2}\begin{array}{c|c|c} 
\mathbf{0} &  I_{k-2} & \mathbf{0} \\
\hline 
0 & \mathbf{0} & 0 \\
\hline
1 &  \mathbf{0} &  0 \end{array}} \right ]=N_3N_4,
\label{eq:jk0_k_odd}
\end{equation} when $k$ is odd and \begin{equation}
J_k(0) = \left [ {\renewcommand{\arraystretch}{1.2}\begin{array}{c|c|c} 
\mathbf{0} & \mathbf{0} & \begin{matrix} 0 & 0 \\ 1 & 1 \end{matrix}\\
\hline 
\begin{matrix} 1 & -1\\ 0 & 1 \end{matrix} &  \mathbf{0} &  \mathbf{0} \\
\hline
\mathbf{0} & I_{k-4} & \mathbf{0} \end{array}} \right ]
(N_3 + E_{(1,3)}) = N_5N_6,
\label{eq:jk0_k_even}
\end{equation}  when $k$ is even. For \eqref{eq:jk0_k_odd} let \begin{equation} Q_3=[e_1,e_3,\ldots,e_{k},e_2,e_4,\ldots,e_{k-1}], \label{eq:jk0_k_odd_q1} \end{equation} then it is easy to verify that both factors in 
\[Q_3^{-1} J_k(0) Q_3 = (Q_3^{-1} N_3 Q_3) (Q_3^{-1} N_4 Q_3) \] are of the desired form. In particular note that $(Q_3^{-1} N_4 Q_3)$ has the zero vector as its last row. Now for \eqref{eq:jk0_k_even}, let \begin{equation}Q_4 = [e_1,e_3,\ldots,e_{k-1},e_2,e_4-e_3,e_6-e_5,\ldots,e_k-e_{k-1}], 
\label{eq:jk0_k_even_q2} 
\end{equation} then both factors in
\[Q_4^{-1} J_k(0) Q_4 = (Q_4^{-1} N_5 Q_4) (Q_4^{-1} N_6 Q_4) \] are of the desired form. In particular $(Q_4^{-1} N_6 Q_4)$ has as its last row the vector $e_1^T$. 

Now suppose the last diagonal block of $N$ is $J_2(0)$. Suppose first that we have $N = \text{Dg}[J,J_k(0),J_2(0)]$ where $J$ is some nilpotent matrix in Jordan canonical form. If $k=1$ the result is immediate by the factorization \eqref{eq:nilpotent_jkj2_1}. If $k \geq 2$ and $k$ is even the result is also immediate by the factorization \begin{equation} \text{Dg}[J_k(0),J_2(0)] = \left [ {\renewcommand{\arraystretch}{1.2}\begin{array}{c|c} 
\mathbf{0} & J_k(0) \\ 
\hline
\begin{matrix} 0 & 0 \\ 0 & 1 \end{matrix} &  \mathbf{0} \end{array}} \right ] 
\left [ {\renewcommand{\arraystretch}{1.2}\begin{array}{cc|c} 
\mathbf{0} & \mathbf{0} & J_2(0) \\ 
\hline
 I_{k-1} & \mathbf{0}  & \mathbf{0} \\ 
\mathbf{0} & 0 &  0 \end{array}} \right ]=N_7N_8. \label{eq:nilpotent_jkj2_even} \end{equation} Suppose $k \geq 2$ and $k$ is odd: we can make use of the factorization as defined in \eqref{eq:nilpotent_jkj2_odd}. By applying the change-of-basis matrix $Q_1$ as defined in \eqref{eq:nilpotent_jkj2_odd_q}, both factors in 
\[Q_1^{-1} \text{Dg}[J_k(0),J_2(0)] Q_1 = (Q_1^{-1} N_1 Q_1) (Q_1^{-1} N_2 Q_1) \] are for the desired form. In particular note that the last row of $Q_1^{-1} N_2 Q_1$ is the zero vector.

Finally suppose $N = \text{Dg}[J,J_2(0),J_2(0),J_2(0)]$. We can make use of the factorization  \begin{equation} \text{Dg}[J_2(0),J_2(0),J_2(0)] = \left [ {\renewcommand{\arraystretch}{1.2}\begin{array}{c|c|c} 
\mathbf{0} &  \mathbf{0} & \begin{matrix} 0 & 0 \\ 0 & 1 \end{matrix} \\ 
\hline
J_2(0) & \mathbf{0} & \mathbf{0} \\
\hline
 \mathbf{0} & J_2(0) & \mathbf{0} \end{array}} \right ]
 \left [ {\renewcommand{\arraystretch}{1.2}\begin{array}{c|c|c} 
\mathbf{0} & \begin{matrix} 1 & 0 \\ 0 & 0 \end{matrix} &  \mathbf{0} \\ 
\hline
\mathbf{0} & \mathbf{0} &  \begin{matrix} 1 & 0 \\ 0 & 0 \end{matrix}  \\
\hline
  J_2(0) & \mathbf{0} & \mathbf{0} \end{array}} \right ] = N_9N_{10}, \label{eq:j2j2j2} \end{equation} and the change-of-basis matrix \begin{equation} Q_5 = [e_1,e_4,e_3,e_6,e_2,e_5], \label{eq:j2j2j2_q} \end{equation} then both factors in 
\[ Q_5^{-1} \text{Dg}[J_2(0),J_2(0),J_2(0)] Q_5 = (Q_5^{-1} N_9 Q_5) (Q_5^{-1} N_{10} Q_5) \] are of the desired form, in particular note that the last row of $Q_5^{-1} N_{10} Q_5$ is the zero vector. 

Now, to prove that the top row of the first factor of $N$ is the zero vector, we can again consider each of the cases above, but this time assuming that the block(s) presented are the first diagonal block(s). Notice that in each of the given factorizations, together with the necessary similarity transformation where needed, the top row of the first factor is the zero vector, which completes the proof. \hfill $\square$

On to the main result of this section, as proved by Laffey. I modify the proof to include the converse, for completeness. I give the proof in full.
\begin{quote}
\begin{thm}[Theorem 1.3 \cite{laffey}]
Let $\mathscr{F}$ be an arbitrary field. The matrix $A \in M_n(\mathscr{F})$ is a product of two nilpotent matrices if and only if it is singular, except when $A \in M_2(\mathscr{F})$ and $A$ is nilpotent and nonzero.
\label{thm:nilpotentProductLaffey}
\end{thm}
\end{quote}
{\bf Proof.} First, suppose $A$ is the product of two nilpotent matrices $N_1,N_2$. Now a nilpotent matrix cannot be full rank, since it has zero as a characteristic value. Furthermore $\text{r}(A) \leq \min(\text{r}(N_1),\text{r}(N_2))$ \cite[Proposition 6.11]{golan}, and it follows that $A$ cannot be full rank, and is therefore singular. 

Now suppose $A$ is singular. If $A$ is nilpotent then the result follows directly by Lemma \ref{lem:j0nilpotent} and \cite[Lemma 3]{wu}. Suppose therefore that $A$ is not nilpotent. By Fitting's lemma \cite[Theorem 5.10]{cullen} we have $A=\text{Dg}[A_0,A_1]$ where $A_0$ is nilpotent and $A_1$ is invertible. Let us consider three mutually exclusive cases in terms of the matrix $A_0$: first suppose $A_0 = J_1(0)$. Let $k=n-1$. Since $A_1$ is similar to $LU$ where $L=(l_{ij})$ is lower triangular and $U=(u_{ij})$ is upper triangular (Theorem 1.1 in \cite{laffey}) we have 
\[A = \begin{bmatrix} 0 & 0 & \cdots & & 0 \\
 l_{11} & 0 & & & \vdots \\
 l_{21} & l_{22} & \ddots & & \vdots \\
 \vdots & & \ddots & 0 & 0 \\
 l_{k1} & \cdots & & l_{kk} & 0
 \end{bmatrix} 
 \begin{bmatrix} 0 &  u_{11} & u_{12} & \cdots & u_{1k} \\
 0 & 0 & u_{21} & & \vdots \\
 \vdots &  & \ddots & \ddots & \vdots \\
 \vdots & &  & 0 & u_{kk} \\
 0 & \cdots & &  & 0
 \end{bmatrix}.\]
 Note that each factor on the right-hand-side has the same rank as $A$, and characteristic polynomial $x^n$ so that it is nilpotent.
 
Now suppose that $A_0$ is similar to $J_2(0)$. Note that $J_2(0)$ is similar to 
\[\begin{bmatrix} 0 & 1 \\ 0 & 0 \end{bmatrix},\] which is the Jordan form preferred by some texts. Let $k=n-2$, and let
\[A_2 = \begin{bmatrix} 
1 	& 0 		& \cdots 	& 		& 0 		& 1\\
0 	& 0 		& \cdots 	& 		& 0 		& 0\\
0 	& l_{11} 	& 0 		& 		& \vdots	& \vdots \\
0	& l_{21} 	& l_{22} 	& \ddots 	& 		&  \\
\vdots & \vdots 	& 		& \ddots 	& 0 		& 0 \\
-1 	& l_{k1} 	& \cdots 	& 		& l_{kk} 	& -1
 \end{bmatrix}
 \begin{bmatrix} 
0 	& 1 		& 0	 	& \cdots		&  		& 0\\
0 	& 0 		&  u_{11} 	& u_{12} 	& \cdots 	& u_{1k}\\
0 	& 0 		& 0 		& u_{21} 	& 		& \vdots \\
\vdots	& \vdots 	&  		& \ddots 	& \ddots 	& \vdots\\
 &  	& 		&  		& 0 		& u_{kk} \\
0 	& 0 		& \cdots 	& 		& 	0 	& 0
 \end{bmatrix} =
 \left [ {\renewcommand{\arraystretch}{1.2}\begin{array}{c|c} 
\begin{matrix} 0 & 1 \\ 0 & 0 \end{matrix} &  \mathbf{0} \\ 
\hline
-E_{(k,2)} & A_1 \end{array}} \right ].
 \]
By Roth's theorem $A_2$ is similar to $A$ \cite{roth}. Explicitly, 
\[ \begin{bmatrix} I & \mathbf{0} \\ -X & I\end{bmatrix}  \left [ {\renewcommand{\arraystretch}{1.2}\begin{array}{c|c} 
\begin{matrix} 0 & 1 \\ 0 & 0 \end{matrix} &  \mathbf{0} \\ 
\hline
-E_{(k,2)} & A_1 \end{array}} \right ] \begin{bmatrix} I & \mathbf{0} \\ X & I\end{bmatrix} =  \left [ {\renewcommand{\arraystretch}{1.2}\begin{array}{c|c} 
\begin{matrix} 0 & 1 \\ 0 & 0 \end{matrix} &  \mathbf{0} \\ 
\hline
\mathbf{0} & A_1 \end{array}} \right ]\] only if there exists a solution to 
\[-E_{(k,2)} = X \begin{bmatrix} 0 & 1 \\ 0 & 0 \end{bmatrix} - A_1X.\] 
Now since $A_1$ is full rank there exists a vector $x \in \mathscr{F}^k$ such that $A_1x = e_k$. It follows that $X = [0,x]$ is a matrix that will satisfy the requirements, proving that the given factorization is valid. 

Let the first factor in the factorization above be $N_1$ and the second $N_2$. Since $N_2$ is upper triangular with only 0 on the diagonal, it is immediately apparent that it is nilpotent. Now consider $N_1$: Let $P=[e_1,e_1-e_n,e_2,e_3,\ldots,e_{n-1}]$, then
\[ P^{-1} N_1 P = \left [ {\renewcommand{\arraystretch}{1.2}\begin{array}{c|c} 
\begin{matrix} 0 & 0 \\ 1 & 0 \end{matrix} &  \begin{matrix} l_{k1} & l_{k2} & \cdots & \qquad & l_{kk} \\ 
-l_{k1} & -l_{k2} & \cdots & & -l_{kk} \end{matrix} \\ 
\hline
\mathbf{0} & \begin{matrix} 
0 & 0 & \cdots & & 0 \\
l_{11} & 0 &  & & 0 \\ 
l_{21} & l_{22} & \ddots & & \vdots \\
\vdots & & \ddots & 0 & 0 \\
l_{1(k-1)} & \cdots & & l_{(k-1)(k-1)} & 0
\end{matrix} \end{array}} \right ]. \]
Since the determinant of a block triangular matrix is the product of the determinants of its diagonal blocks \cite[Proposition 11.12]{golan}, the characteristic polynomial of $N_1$ is the product of the characteristic polynomials of its diagonal blocks. It is easy to verify that the characteristic polynomial of $N_1$ is therefore $x^2 \cdot x^k = x^n$, confirming that it is nilpotent, and concluding the proof of the result in this case. 

It remains to prove the result for the case where $A_1$ is of order $k \times k$ where $k \leq n-3$, so that $A_0$ is of order $3 \times 3$ or larger. I will now show that with a slight modification of the factorization given in the original proof the result remains valid. For this purpose I make use of Lemma \ref{lem:col0row0} whereby we can assume $A_0 = N_1N_2$, where $N_1$ is nilpotent and its last column is the zero vector, and $N_2$ is nilpotent and its last row is either the zero vector or $e_1^T$.

Now we have
\[ \begin{bmatrix} A_0 & \mathbf{0} \\ B & A_1\end{bmatrix} =  \left [ {\renewcommand{\arraystretch}{1.2}\begin{array}{c|c} 
N_1 & \mathbf{0} \\
\hline
\begin{matrix} 
0 & \cdots & 0 & l_{11} \\
0 &  & \vdots &  l_{21} \\ 
\vdots &  & & \vdots \\
0 & \cdots & 0 & l_{k1} 
\end{matrix} & \begin{matrix} 
0 & 0 & \cdots & 0 \\
l_{22} & 0 &  &  0 \\ 
\vdots & \ddots & \ddots & \vdots \\
l_{k2} & \cdots & l_{kk} & 0 
\end{matrix}
\end{array}} \right ]
\left [ {\renewcommand{\arraystretch}{1.2}\begin{array}{c|c} 
N_2 &\begin{matrix} 
0 & \cdots &   & 0  \\
\vdots &  &    & \vdots \\ 
0 & \cdots  &  & 0 \\
u_{11} & u_{12} &  \cdots & u_{1k} 
\end{matrix}  \\
\hline
\mathbf{0} &
 \begin{matrix}
0  & u_{22} &  \cdots & u_{2k} \\
\vdots & \ddots & \ddots    &  \vdots \\ 
 \quad \, &  & 0 & u_{kk} \\
0 & \cdots & 0 & 0 
\end{matrix}
\end{array}} \right ].\]
If the last row of $N_2$ is zero then $B=\mathbf{0}$ and it follows that $A$ is similar to the given factorization, which is the desired result. Suppose the last row of $N_2$ is $e_1^T$, then   
\[ B= \begin{bmatrix} l_{11} &0 & \cdots & 0 \\
l_{21} &0 & \cdots & 0 \\
\vdots  & \vdots & & \vdots \\
l_{k1}&0 & \cdots & 0 \end{bmatrix}. \] Now we may use Roth's theorem to prove the result \cite{roth}. Explicitly, 
\[ \begin{bmatrix} I & \mathbf{0} \\ -X & I\end{bmatrix} \begin{bmatrix} A_0 & \mathbf{0} \\ B & A_1\end{bmatrix} \begin{bmatrix} I & \mathbf{0} \\ X & I\end{bmatrix} = \begin{bmatrix} A_0 & \mathbf{0} \\ \mathbf{0} & A_1\end{bmatrix}\] only if there exists a solution to $B = XA_0 - A_1X$. Let $X$ be the matrix with $-1/u_{11}$ in entry $(1,1)$ and zero elsewhere. Then $XA_0 = 0$ since the first row of $A_0$ is zero (by Lemma \ref{lem:col0row0} we may assume the first row of $N_1$ is the zero vector), and $-A_1X = B$, which yields the desired result. This proves that $A$ is similar to the given factorization.

It remains to prove that the factors are nilpotent. But now, as mentioned before, since the determinant of a block triangular matrix is the product of the determinants of its diagonal blocks \cite[Proposition 11.12]{golan}, the characteristic polynomial of each factor is the product of the characteristic polynomials of its diagonal blocks. Notice that for both these factors the diagonal blocks are all nilpotent, and therefore both factors are nilpotent, which concludes the proof. \hfill $\square$

\end{document}